\documentclass[11pt]{article}
\usepackage{amsmath}
\usepackage{amssymb}
\usepackage{amsthm,amsxtra}
\usepackage{epsf}
\usepackage{eepic}
\newlength{\mytopmargin}
\newlength{\myleftmargin}
\newtheorem{thm}{Theorem}
\newtheorem{cor}{Corollary}
\newtheorem{lemma}{Lemma}
\newtheorem{prop}{Proposition}
\begin{document}
\vspace{4cm}
\noindent
{\bf Jacobians and rank 1 perturbations
relating to unitary Hessenberg matrices}

\vspace{5mm}
\noindent
Peter J.~Forrester${}^*$\footnote{Supported by the Australian Research
Council}
 and Eric M.~Rains${}^\dagger$

\noindent
${}^*$Department of Mathematics and Statistics,
University of Melbourne, \\
Victoria 3010, Australia ;
${}^\dagger$
Department of Mathematics, University of California, Davis, CA 95616, USA

\small
\begin{quote}
In a recent work Killip and Nenciu gave random recurrences for the
characteristic polynomials of certain unitary and real orthogonal upper
Hessenberg matrices. The corresponding eigenvalue p.d.f.'s are
$\beta$-generalizations of the classical groups. Left open was the direct
calculation of certain Jacobians. We provide the sought direct calculation.
Furthermore, we show how a multiplicative rank 1 perturbation of the
unitary Hessenberg matrices provides a joint eigenvalue p.d.f.~generalizing
the circular $\beta$-ensemble, and we show how this joint density is related
to known inter-relations between circular ensembles. Projecting the joint
density onto the real line leads to the derivation of a random three-term
recurrence for polynomials with zeros distributed according to the
circular Jacobi $\beta$-ensemble.
\end{quote}

\section{Introduction}

Consider the classical group $U(N)$ of $N \times N$ unitary matrices. There is
a unique measure $d_HU$ --- the Haar measure --- which is invariant under
both left and right multiplication by a fixed unitary matrix, thus
giving a uniform distribution on the group.
The corresponding eigenvalue probability
density function (p.d.f.) has the explicit form (see e.g.~\cite{Fo02})
\begin{equation}\label{1.1}
{1 \over (2 \pi)^N N!} \prod_{1 \le j < k \le N} | e^{i \theta_k} -
e^{i \theta_j} |^2,
\end{equation}
and this in turn is of fundamental importance in recent applications of
random matrix theory to combinatorial models \cite{Ra98,BR01a}, analytic
number theory \cite{KS00a} and the quantum many body problem
\cite{FFGW02}.

A basic question is how to best sample from (\ref{1.1}). Until very
recently, the only method available has been to first generate a member
of $U(N)$ according to the Haar measure, by for example applying the
Gram-Schmidt orthogonalization procedure to the columns of an
$N \times N$ complex Gaussian matrix, then to calculate the
eigenvalues of the resulting matrix. However, inspired by recent work of
Dumitriu and Edelman \cite{DE02}, this situation has been dramatically
improved upon by Killip and Nenciu \cite{KN04}. Thus augmenting ideas
from \cite{DE02} with results from the theory of orthogonal polynomials
on the unit circle, these authors have provided an explicit unitary
Hessenberg matrix, with positive elements on the subdiagonal, which
has for its eigenvalue p.d.f.~the $\beta$-generation of (\ref{1.1}),
\begin{equation}\label{1.2}
{1 \over C_{\beta N} } \prod_{1 \le j < k \le N} | e^{i \theta_k} -
e^{i \theta_j} |^\beta, \qquad \qquad
C_{\beta N} = (2 \pi)^N {\Gamma (\beta N/2 + 1) \over
(\Gamma(\beta/2))^N }.
\end{equation}
In general the characteristic polynomial 
$\chi_N(\lambda)$ of such matrices
can be calculated from the coupled recurrences
\begin{eqnarray}
\chi_k(\lambda) & = & \lambda \chi_{k-1}(\lambda) - \bar{\alpha}_{k-1}
\tilde{\chi}_{k-1}(\lambda) \nonumber \\
\tilde{\chi}_k(\lambda) & = & \tilde{\chi}_{k-1}(\lambda) - \lambda
\alpha_{k-1} \chi_{k-1}(\lambda) \label{ck}
\end{eqnarray}
$(k=1,\dots,N)$ where $\chi_0(\lambda) =  \tilde{\chi}_{0}(\lambda) = 1$
and furthermore 
$\tilde{\chi}_k(\lambda) = \lambda^k \bar{\chi}_k(1/\lambda)$.
For the unitary Hessenberg matrix relating to (\ref{1.2}),
the parameters $\{ \alpha_j\}_{j=0,\dots,N-1}$ are random variables with
distributions specified in \cite{KN04} (see (\ref{he}) below). As a
consequence of this result the joint distribution (\ref{1.1}), or
more generally (\ref{1.2}), can be sampled by simply iterating (\ref{ck})
to generate $\chi_N(\lambda)$, then computing its roots.

The problem of efficiently sampling from the p.d.f.
\begin{equation}\label{1.3}
{1 \over C_N(a,b;\beta) } \prod_{l=1}^N |1 - e^{i \theta_l} |^{2a+1}
|1 + e^{i \theta_l} |^{2b+1} 
\prod_{1 \le j < k \le N} |e^{i \theta_j} - e^{i \theta_k} |^\beta
|1 - e^{i (\theta_j + \theta_k)} |^\beta, \quad
(0 \le \theta_l \le \pi)
\end{equation}
was solved according to the same strategy in \cite{KN04}. Here the
underlying unitary Hessenberg matrix is real orthogonal with determinant
$+1$, and thus the characteristic polynomial $\chi_{2N}(\lambda)$
has real coefficients.
In this
case $ \tilde{\chi}_k(\lambda) = \lambda^k \chi_k(1/\lambda)$ and so only the
first of the recurrences in (\ref{ck}) 
is required. 
Note that the eigenvalues of a real orthogonal matrix with determinant
$+1$ come in complex conjugate pairs $e^{ \pm i \theta}$;
(\ref{1.3}) is the joint
distribution of those with $0 < \theta < \pi$.
The case $\beta = 2$, $(a,b)=(\pm {1 \over 2},
\pm {1 \over 2})$ (the signs chosen appropriately) 
of (\ref{1.3})
gives the
eigenvalue p.d.f.~ for matrices from the real orthogonal and symplectic
classical groups with Haar measure (see e.g.~\cite{Fo02}). Like their
counterparts from $U(N)$, such random matrices are of fundamental importance
in applications of random matrix theory to combinatorial models \cite{Ra98,
BR01a},
analytic number theory \cite{KS00b} and the quantum many body problem
\cite{FFG03}.

In the work \cite{KN04}, Killip and Nenciu left open two questions concerning 
the direct computation of certain Jacobians, one relating to unitary
Hessenberg matrices corresponding to (\ref{1.2}), and the other to real
orthogonal Hessenberg matrices corresponding to (\ref{1.3}). Earlier,
Dumitriu and Edelman \cite{DE02} had left open an analogous question in the
case of tridiagonal matrices corresponding to the Gaussian $\beta$-ensemble
p.d.f.
$$
{1 \over G_{\beta N} } \prod_{l=1}^N e^{-x_l^2/2} \prod_{1 \le j < k \le N}
| x_k - x_j |^\beta.
$$
It was remarked in \cite{DE02} that one of the present authors (PJF) had
communicated a direct derivation of the sought Jacobian for the change of
variables from the elements of a tridiagonal matrix, to the eigenvalues
and the first component of the eigenvectors. A primary purpose of this
article is to show how a similar approach can be used to answer the two
questions left open in \cite{KN04}. We begin in Section 2 by presenting
the calculation for the Jacobian in the case of a tridiagonal matrix.
In Section 3 this calculation is extended to provide a direct calculation
of Jacobians relating to unitary and real orthogonal Hessenberg matrices.
Also shown is how portions of the working in \cite{KN04} reliant on the
theory of orthogonal polynomials on the unit circle, can alternatively
be derived within a matrix setting. In Section 4, it
is shown how a certain rank 1 multiplicative
perturbation of unitary matrices leads to the
derivation of a joint eigenvalue p.d.f.~generalizing (\ref{1.2}). An
integration formula associated with this p.d.f.~is discussed, which
in turn is shown to include as special cases known
inter-relations between circular ensembles. Furthermore, the
multiplicative perturbation is used to give an alternative derivation of
these inter-relations.

The Cayley transformation of the distributions obtained in Section 4, projecting the
unit circle to the real line, are studied in Section 5. This leads to a
random three-term recurrence for the projection  onto the
real line of polynomials with zeros distributed according to
\begin{equation}\label{Mn}
{1 \over M_N(a;c) } \prod_{l=1}^N |1 - e^{i \theta_l} |^a
\prod_{1 \le j < k \le N} |e^{i \theta_k} - e^{i \theta_j} |^{2c},
\end{equation}
which with $2c = \beta$ is known as the circular Jacobi $\beta$-ensemble \cite{Fo02}.
In the case $a=0$ this recurrence scheme is distinct from the
scheme (\ref{ck}).

\section{Calculation of a Jacobian for tridiagonal matrices}
\setcounter{equation}{0}
Let
\begin{equation}\label{T}
T = \left [ \begin{array}{ccccc} a_n & b_{n-1} & & & \\
b_{n-1} & a_{n-1} & b_{n-2} & & \\
 & b_{n-2} & a_{n-2} & b_{n-3} & \\
& \ddots & \ddots & \ddots & \\
&&b_{2} & a_2 & b_{1} \\
& & & b_{1} & a_1 \end{array} \right ]
\end{equation}
be a general real symmetric tridiagonal matrix. The problem posed in
\cite{DE02} 
is to compute the Jacobian for the change of variables from the description
of $T$ in terms of its entries, to the description in terms of eigenvalues
and variables relating to its eigenvectors.

As is well known, and easy to see by direct substitution,
for each eigenvalue $\lambda_k$ and corresponding
eigenvector $\vec{v}_k$, once 
the 1st component $\vec{v}_k^{(1)} =: q_k$, $q_k >0$,
of $\vec{v}_k$ is specified,
the other components are then fully determined by $\{\lambda_k\}$ and the
elements of $T$. 
However only $n-1$ of these components are independent due to the relation
\begin{equation}\label{qq}
\sum_{k=1}^n q_k^2 = 1,
\end{equation}
which itself is a consequence of $T$ being symmetric and thus orthogonally
diagonalizable.
Thus the $2n-1$ variables 
\begin{equation}\label{1.62}
\vec{a} := (a_n,a_{n-1},\dots,a_1), \qquad
\vec{b} := (b_{n-1},\dots, \vec{b}_1)
\end{equation}
can be put into 1-to-1 correspondence with the $2n-1$ variables
\begin{equation}\label{1.63}
\vec{\lambda} := (\lambda_1,\dots, \lambda_n), \qquad
\vec{q} := (q_1,\dots, q_{n-1})
\end{equation}
where $\lambda_1 > \cdots > \lambda_n$ and $q_i > 0$. The Jacobian for the
change of variables from (\ref{1.62}) to (\ref{1.63}) can be computed
directly using the method of wedge products (for an introduction to the
use of this technique in random matrix theory see \cite{Fo02}).

We will first isolate results required in the course of the calculation.

\begin{prop}\label{p3}
Let $(X)_{ij}$ denote the $ij$ entry of the matrix $X$. We have
\begin{equation}\label{a}
((T - \lambda I_n)^{-1})_{11} = \sum_{j=1}^n {q_j^2 \over \lambda_j - 
\lambda}.
\end{equation}
Also
\begin{equation}\label{b}
\prod_{1 \le i < j \le n}(\lambda_i - \lambda_j)^2
= {\prod_{i=1}^{n-1} b_i^{2i} \over \prod_{i=1}^n q_i^2}
\end{equation}
and
\begin{equation}\label{c}
\det \Big [ [\lambda_k^j - \lambda_n^j]_{j=1,\dots,2n-1 \atop
k = 1,\dots, n-1}
[j \lambda_k^{j-1}]_{j=1,\dots,2n-1 \atop k=1,\dots,n} \Big ]
= \prod_{1 \le j < k \le n} (\lambda_k - \lambda_j)^4.
\end{equation}
\end{prop}

\noindent
Proof. \quad The identity (\ref{a}), which is well known, follows by writing
the matrix entry as an inner product, and decomposing the vectors in this
inner product as eigenvectors. The identity (\ref{b}) is contained in
\cite{DE02}. It can be derived from (\ref{a}) by using the fact that for a
general $n \times n$ non-singular matrix
\begin{equation}\label{X}
(X^{-1})_{11} = {\det X_{n-1} \over \det  X},
\end{equation}
where $X_{n-1}$ denotes the bottom right $n-1 \times n-1$ submatrix of $X$,
introducing the corresponding characteristic polynomials
$P_{n-1}(\lambda)$, $P_n(\lambda)$, and making use of the three term
recurrence
$$
P_k(\lambda) = (\lambda - a_k) P_{k-1}(\lambda) - b_{k-1}^2 P_{k-2}(\lambda),
\quad P_0(\lambda) := 1.
$$
For the identity (\ref{c}), note that both sides are homogeneous symmetric
polynomials of degree $2n(n-1)$. Furthermore, the determinant and its first
three derivatives with respect to $\lambda_1$ vanish at $\lambda_1 =
\lambda_2$. As a consequence, it follows that the determinant must in fact
be proportional to the fourth power of the product of differences as given
in the r.h.s. The fact that the proportionality constant is unity
follows by comparing coefficients of $(\lambda_1^0 \lambda_2^1 \cdots
\lambda_n^{n-1})^4$ on both sides.
\hfill   $\square$

\medskip
The Jacobian can now be computed according to the following result.

\begin{thm}\label{th1}
The Jacobian for the change of variables (\ref{1.62}) to (\ref{1.63}) is
equal to
\begin{equation}\label{qb}
{1 \over q_n} {\prod_{i=1}^{n-1} b_i \over \prod_{i=1}^n q_i}.
\end{equation}
\end{thm}

\noindent
Proof. \quad Rewriting (\ref{a}) in the form
\begin{equation}\label{1.70'}
((I_n - \lambda T )^{-1})_{11} = \sum_{j=1}^n
{q_j^2 \over 1 - \lambda \lambda_j },
\end{equation}
recalling the explicit form of $T$ from (\ref{T}), and equating successive
powers of $\lambda$ on both sides gives
\begin{eqnarray}\label{2.q}
1 & = &  \sum_{j=1}^n q_j^2 \nonumber \\
a_n  & = & \sum_{j=1}^n q_j^2 \lambda_j \nonumber \\ 
{*} + b_{n-1}^2 &  = &  \sum_{j=1}^n q_j^2 \lambda_j^2 \nonumber \\
{*} + a_{n-1} b_{n-1}^2   & = & \sum_{j=1}^n q_j^2 \lambda_j^3 \nonumber \\
{*} + b_{n-2}^2 b_{n-1}^2 &  = &  \sum_{j=1}^n q_j^2 \lambda_j^4 \nonumber \\
{*} + a_{n-2} b_{n-2}^2 b_{n-1}^2  &  = & \sum_{j=1}^n q_j^2 \lambda_j^5 
\nonumber\\
\vdots &   & \vdots \nonumber \\
{*} + a_{1}b_1^2 \cdots  b_{n-2}^2 b_{n-1}^2   & = &
\sum_{j=1}^n q_j^2 \lambda_j^{2n-1}.
\end{eqnarray}
Here the $*$ denotes terms involving only variables already having
appeared on the l.h.s.~of preceding equations. 
Thus the variables
$a_n, b_{n-1}, a_{n-1}, b_{n-2}, \dots$ occur in a triangular
structure. Upon taking differentials, the first of these equations
implies
$$
q_n dq_n = - \sum_{j=1}^{n-1} q_j dq_j.
$$
For the differentials of the remaining equations, we use this to substitute
for $dq_n$, and then take wedge products of both sides. On the l.h.s., the
triangular structure gives
\begin{equation}\label{71}
2^{n-1} \prod_{j=1}^{n-1}b_j^{4j - 3} d\vec{a} \wedge d\vec{b}
\end{equation}
where
$$
d\vec{a} := \bigwedge_{j=1}^n da_j, \qquad d\vec{b} 
:= \bigwedge_{j=1}^{n-1} db_j.
$$
On the r.h.s.~the wedge product operation yields
\begin{equation}\label{72}
2^{n-1} q_n^2 \prod_{j=1}^{n-1} q_j^3
\det \Big [ [\lambda_k^j - \lambda_n^j]_{j=1,\dots,2n-1 \atop
k = 1,\dots, n-1}
[j \lambda_k^{j-1}]_{j=1,\dots,2n-1 \atop k=1,\dots,n} \Big ]
d \vec{\lambda} \wedge d \vec{q}
\end{equation}
where
$$
d \vec{\lambda} := \bigwedge_{j=1}^nd\lambda_j, \qquad
d\vec{q} := \bigwedge_{j=1}^{n-1} dq_j
$$
(a common factor $2q_k$ has been removed from column $k$,
$k=1,\dots,n-1$ of the determinant, as has a common factor
$q_k^2$ from columns $n-1+k$, $k=1,\dots,n$).

By definition the Jacobian $J$ satisfies
\begin{equation}\label{2.14'}
d \vec{a} \wedge d \vec{b} = J d \vec{\lambda} \wedge d \vec{q}.
\end{equation}
Equating (\ref{71}) and (\ref{72}), and using (\ref{c}) to evaluate
the determinant then shows
$$
J = {1 \over q_n}
{\prod_{j=1}^{n-1} b_j \over \prod_{j=1}^n q_j}
\Big ( {\prod_{j=1}^n q_j^2 \over \prod_{j=1}^{n-1}
b_j^{2j-1}} \Big )^2 \prod_{1 \le j < k \le N}
(\lambda_k - \lambda_j)^4.
$$
Recalling (\ref{b}) gives the form of $J$ (\ref{qb}).
\hfill $\square$

\medskip
In \cite{DE02} indirect methods are used to derive (\ref{qb}) but
with the factor $1/q_n$ not present. As noted in
\cite{KN04}, the reasoning of \cite{DE02} is most suited to working
with the variables $\mu_j = q_j^2$, and doing so eliminates this apparent
discrepancy.

\section{Calculation of a Jacobian for unitary and real orthogonal
Hessenberg matrices}

\setcounter{equation}{0}

\subsection{Preliminaries}
In general a unitary upper triangular Hessenberg matrix $H = 
[H_{i,j}]_{i,j=
0,\dots,n-1}$ with positive elements along the sub-diagonal is
parametrized by $n-1$ complex numbers $\alpha_0,\alpha_1,\dots,
\alpha_{n-2}$ with $|\alpha_j| = 1$ and a further complex number
$\alpha_{n-1}$ with $|\alpha_{n-1}| < 1$. Setting $\alpha_{-1} :=
-1$, $\rho_j := \sqrt{1 - |\alpha_j|^2}$ $(j=0,\dots,n-2)$, one
can check that if the diagonal entries are specified as $H_{i,i} = -
\alpha_{i-1} \bar{\alpha}_i$, and subdiagonal entries as
$H_{i+1,i} = \rho_i$, then the remaining non-zero entries are
given by
\begin{equation}\label{Hr}
H_{i,j} = - \alpha_{i-1} \bar{\alpha}_j \prod_{l=i}^{j-1} \rho_l, \qquad
i < j.
\end{equation}

Let $\lambda_j = e^{i \theta_j}$ ($j=1,\dots,n$) denote the eigenvalues
of $H$ and let $q_j$ denote the modulus of the first component of the
corresponding normalized eigenvectors (the $\{q_j\}$ thus satisfy
(\ref{qq})). With the $\{\theta_j\}$ ordered, there is an
invertible 1-to-1 correspondence with the parameters
$\{\alpha_j\}_{j=0,\dots,n-1}$. Our interest is to directly compute
the Jacobian for the change of variables from $\{\alpha_j\}_{j=0,\dots,n-1}$
to $\{\theta_i\}_{i=1,\dots,n}$, $\{q_i\}_{i=1,\dots,n-1}$.

In preparation for the derivation of this result, note that with
$H_k$ denoting the top $k \times k$ block of $H$,
$$
\chi_k(\lambda) := \det(\lambda I_k - H_k)
$$
satisfies (\ref{ck}) (see e.g.~\cite{Gr93}). 
Also of interest is a variant of the
characteristic polynomial of the bottom $k \times k$ submatrix. In relation
to this, note that with the involution $\alpha_j \mapsto -
\bar{\alpha}_j \alpha_{n-1}$ $(j=0,\dots,n-2)$, the bottom $k \times k$
submatrix, after reflection in the anti-diagonal, becomes equal to the
top $k \times k$ submatrix but with $\alpha_j \mapsto \alpha_{n-2-j}$
($j=0,\dots,n-2$). Let the characteristic polynomial of
the bottom $k \times k$ submatrix with the replacements
$\alpha_j \mapsto -
\bar{\alpha}_j \alpha_{n-1}$ $(j=0,\dots,n-2)$ be denoted
$\chi_k^b(\lambda)$. We see that this polynomial satisfies the recurrence
(\ref{ck}) with 
$\alpha_j$ replaced by 
$- \bar{\alpha}_{n-2-j} \alpha_{n-1}$ in (\ref{ck}),
\begin{eqnarray}
\chi_k^b(\lambda) & = & \lambda \chi_{k-1}^b(\lambda) + 
\alpha_{n-1-k} \bar{\alpha}_{n-1}
\tilde{\chi}_{k-1}^b(\lambda) \nonumber \\
\tilde{\chi}_k^b(\lambda) & = & \tilde{\chi}_{k-1}^b(\lambda) + \lambda
\bar{\alpha}_{n-1-k}\alpha_{n-1} \chi_{k-1}^b(\lambda) \label{c.2}
\end{eqnarray}
$(k=1,\dots,n)$ where $\chi^b_0(\lambda) = \tilde{\chi}^b_0(\lambda) =1$ and
$\tilde{\chi}_k^b(\lambda) = \lambda^k \bar{\chi}_k^b(1/\lambda)$. These
recurrences can be used to derive the analogue of (\ref{b})
\cite{KN04}

\begin{prop}\label{H3}
For the unitary Hessenberg matrix specified by (\ref{Hr}) and surrounding
text, we have
\begin{equation}\label{3H}
\prod_{1 \le i < j \le n} |\lambda_i - \lambda_j |^2 =
{\prod_{l=0}^{n-2} (1 - |\alpha_l|^2)^{n-1-l} \over
\prod_{j=1}^n q_j^2 }.
\end{equation}
\end{prop}

\noindent
Proof. \quad We follow the strategy sketched to prove (\ref{b}) in the
proof of Proposition \ref{p3}. Analogous to (\ref{a}) we have
\begin{equation}\label{aH}
((H- \lambda I_n)^{-1})_{11} = \sum_{j=1}^n {q_j^2 \over \lambda_j
- \lambda}
\end{equation}
where $q_j$ and $\lambda_j$ relate to $H$ as specified below
(\ref{Hr}).
Using (\ref{X}), (\ref{aH}) can  be rewritten as
\begin{equation}\label{3.18'}
{\chi_{n-1}^b(\lambda)  |_{\alpha_j \mapsto - \bar{\alpha}_j
\alpha_{n-1}} \over \prod_{i=1}^n( \lambda - \lambda_i) } =
\sum_{j=1}^n {q_j^2 \over \lambda - \lambda_j},
\end{equation}
which in turn implies
$$
\prod_{i=1}^n q_i^2 = {\prod_{i=1}^n| \chi_{n-1}^b(\lambda_i)| 
|_{\alpha_j \mapsto - \bar{\alpha}_j
\alpha_{n-1}}
\over
\prod_{1 \le i < j \le n} |\lambda_i - \lambda_j|^2 }.
$$
From this we see (\ref{3H}) follows if we can show
\begin{equation}\label{g.ka}
\prod_{i=1}^n| \chi_{n-1}^b(\lambda_i)| = \prod_{l=0}^{n-2} (1 -
| \alpha_l|^2 )^{n-1-l}.
\end{equation}

To establish (\ref{g.ka})
we will use (\ref{c.2}). With $\lambda_j^{(p)}$ denoting
the $j$th zero of $\chi_p^b(\lambda)$, it follows from (\ref{c.2}) that
\begin{eqnarray}\label{g.k}
\chi_k^b(1/\bar{\lambda}_j^{(k)}) & = & {1 \over \bar{\lambda}_j^{(k)}}
(1 - |\alpha_{n-k-1}|^2) \chi_{k-1}^b(1/\bar{\lambda}_j^{(k)}) \nonumber \\
\tilde{\chi}_k^b(\lambda_j^{(k)}) & = & (1 - |\alpha_{n-k-1}|^2)
\tilde{\chi}_{k-1}^b(\lambda_j^{(k)}) .
\end{eqnarray}
Introducing the factorizations
$$
\chi_{k-1}^b(x) = \prod_{i=1}^{k-1} (x - \lambda_i^{(k-1)}), \qquad
\tilde{\chi}_{k}^b(x)= \prod_{i=1}^k (1 - x \bar{\lambda}_i^{(k)})
$$
we deduce from (\ref{g.k}) that
\begin{eqnarray*}
\prod_{i=1}^k \chi_k^b(1/\bar{\lambda}_i^{(k)}) & = &
(1 - |\alpha_{n-k-1}|^2)^k \prod_{i=1}^k(1/\bar{\lambda}_i^{(k)})^k
\prod_{j=1}^{k-1}\tilde{\chi}_{k-1}^b(\lambda_j^{(k-1)}) \\
\prod_{i=1}^k \tilde{\chi}_k^b({\lambda}_i^{(k)}) & = &
(1 - |\alpha_{n-k-1}|^2)^k \prod_{j=1}^{k-1}(\bar{\lambda}_j^{(k-1)})^{k-1}
{\chi}_{k-1}^b(1/\bar{\lambda}_j^{(k-1)}).
\end{eqnarray*}
These latter two equations together imply
$$
\prod_{i=1}^k (\bar{\lambda}_i^{(k)})^k 
{\chi}_{k}^b(1/\bar{\lambda}_i^{(k)}) = \prod_{l=0}^{k-1}(1 - 
|\alpha_{n-l}|^2)^{l+1}.
$$
Making further use of the first equation in (\ref{g.k}), setting $k=n$,
and noting $|\lambda_i^{(n)}|=1$ gives (\ref{g.ka}).
\hfill $\square$ 

\medskip
In \cite{KN04} (\ref{3H}) is derived using a different strategy relating
to the determinant of the Toeplitz matrix formed from the moments of
the underlying measure.

Also required is a determinant evaluation analogous to (\ref{c}).

\begin{prop}
We have
\begin{equation}\label{TN}
\det \left [ \begin{array}{cc}
\left [ \begin{array}{c}\lambda_k^j - \lambda_n^j  \\ 
\lambda_k^{-j} - \lambda_n^{-j} \end{array}
\right ]_{j,k=1,\dots,n-1} 
&
\left [ \begin{array}{c}j \lambda_k^j \\ 
-j \lambda_k^{-j} 
\end{array}
\right ]_{j=1,\dots,n-1 \atop k=1,\dots,n} 
\\
{}[ \lambda_k^n - \lambda_n^n]_{k=1,\dots,n-1} &
[n \lambda_k^n ]_{k=1,\dots,n} \end{array} \right ]
= {\prod_{1 \le j < k \le n} (\lambda_k - \lambda_j)^4 \over
\prod_{l=1}^n \lambda_l^{2n-3} }.
\end{equation}
\end{prop}

\noindent
Proof. \quad By inspection the determinant is a symmetric function of
$\lambda_1,\dots,\lambda_n$ which is homogeneous of degree $n$. Upon
multiplying columns 1 and columns $n$ by $\lambda_1^{2n-3}$ we see that
the determinant becomes a polynomial in $\lambda_1$, so it must be of
the form
$$
{1 \over \prod_{l=1}^n \lambda_l^{2n-3} } p(\lambda_1,\dots,\lambda_n)
$$
where $p$ is a symmetric polynomial of $\lambda_1,\dots,\lambda_n$ of
degree $2n(n-1)$.

We see immediately that the determinant vanishes when $\lambda_1 =
\lambda_2$. Furthermore, it is straightforward to verify that its
derivatives $(\lambda_1 {\partial \over \partial \lambda_1} )^j$
$(j=1,2,3)$ also vanish when $\lambda_1 = \lambda_2$. The polynomial
$p$ must thus contain as a factor
$\prod_{1 \le j < k \le n} (\lambda_k - \lambda_j)^4$. As this is of degree
$2n(n-1)$, it follows that the determinant must in fact be proportional
to (\ref{TN}).

On the r.h.s.~of (\ref{TN}), the coefficient of 
$\prod_{l=1}^n \lambda_l^{4(l-1)-2n+3}$ is unity. In the determinant, let us 
add $(n-1)$ times the first column to the $n$th column. Then we see that
the coefficient of $\lambda_1^{-2n+3}$ is given by the cofactor coming
from multiplying together the $(2n-2,1)$ and $(2n-4,n)$ elements.
In the cofactor we add $(n-2)$ times the first column to the
$(n-1)$st column. The coefficient of $\lambda_1^{-2n+7}$ is given by the
cofactor coming from multiplying together the $(2n-3,1)$ and
$(2n-5,n-1)$ elements. Proceeding in this manner we see that the
coefficient of $\prod_{l=1}^n \lambda_l^{4(l-1)-2n+3}$ is also unity
in the determinant.
\hfill $\square$

\medskip
As remarked in the Introduction, the approach to the ensemble
(\ref{1.3}) given in \cite{KN04} is via $2n \times 2n$  ($n=N$) real orthogonal
Hessenberg matrices with determinant $+1$. The elements being real
implies $\{\alpha_j\}_{j=0,\dots,2n-1}$ are real, while the determinant
equalling $+1$ implies $\alpha_{2n-1} = -1$. Thus there are $2n-1$
independent real parameters $\alpha_0,\dots, \alpha_{2n-2}$. In the 
corresponding eigen-decomposition, there are $n$ independent eigenvalues
$\lambda_j = e^{i \theta_j}$ $(j=1,\dots,n, \: 0 \le \theta_j < \pi)$ and
$n-1$ independent variables $q_j$ $(j=1,\dots,n-1)$ where
${1 \over 2} q_j^2$ is the square of the first component of both the
eigenvalues $\lambda_j$ and $\bar{\lambda}_j$. Left open in
\cite{KN04} is the problem of a direct calculation of the corresponding
Jacobian. For this the analogue of Proposition \ref{H3}  is required.

\begin{prop}\label{I3}
\cite{KN04}
For a $2n \times 2n$ real orthogonal Hessenberg matrix of determinant
$+1$, parametrized in terms of the real parameters
$\{\alpha_i\}_{i=0,\dots,2n-2}$, $|\alpha_i| < 1$, we have
\begin{eqnarray}
&& \prod_{i=1}^n|\lambda_i - {1 \over \lambda_i} |
\prod_{1 \le i < j \le n} |\lambda_i - \lambda_j|^2 |\lambda_i -
1/ \lambda_j |^2 =
2^n {\prod_{l=0}^{2n-2}(1 - \alpha_l^2)^{(2n-1-l)/2} \over
\prod_{i=1}^n q_i^2 }  \label{ii.3} \\
&& \prod_{j=1}^n|1 - \lambda_j|^2 = 2 \prod_{k=0}^{2n-2}(1 - \alpha_k),
\qquad
\prod_{j=1}^n|1 + \lambda_j|^2 = 2 \prod_{k=0}^{2n-2}
(1 + (-1)^k \alpha_k).  \label{iii.3}
\end{eqnarray}
\end{prop}

\noindent
Proof. \quad Denoting the Hessenberg matrix in question by $H$, the analogue
of (\ref{3.18'}) reads
\begin{equation}\label{3.22'}
((I_{2n} - \lambda H)^{-1})_{11} = {1 \over 2}
\sum_{j=1}^n q_j^2 \Big ( {1 \over 1 - \lambda \lambda_j} +
{1 \over 1 - \lambda \bar{\lambda_j}} \Big ).
\end{equation}
Analogous to the reasoning underlying (\ref{3.18'}), the l.h.s.~is equal to
$\chi_{2n-1}^b(1/\lambda)/ \lambda \chi_{2n}^b(1/\lambda)$. We thus have
$$
\Big | {\chi_{2n-1}^b(\lambda_j) \over \chi_{2n}'(\lambda_j)}
\Big | = {1 \over 2} q_j^2 \qquad (j=1,\dots,2n)
$$
where $\lambda_{j+n} = \bar{\lambda}_j$, $q_{j+n} = q_j$. Taking the product
over $j=1,\dots,2n$, making use of (\ref{g.ka}), then taking the square
root gives (\ref{ii.3}). For the results (\ref{iii.3}), one notes
$$
\prod_{j=1}^n|1 - \lambda_j|^2 = \chi_{2n}(1), \qquad
\prod_{j=1}^n|1 + \lambda_j|^2 = \chi_{2n}(-1),
$$
while from (\ref{1.3}) $\chi_{k+1}(\lambda) |_{\lambda = \pm 1} =
(\lambda - \alpha_k \lambda^k) \chi_k(\lambda)  |_{\lambda = \pm 1}$.
\hfill $\square$

\medskip
We remark that in \cite{KN04} (\ref{ii.3}) is deduced by making use of
(\ref{3H}), which in turn is derived using formulas relating to the
underlying measure. Our derivation of (\ref{iii.3}) is the same as that
in \cite{KN04}.

We must make note too of a further determinant evaluation.

\begin{prop}
We have
\begin{eqnarray}\label{Q0}
&&\det \Big [ [ \lambda_k^j + \lambda_k^{-j} - (\lambda_k^j + \lambda_k^{-j})
]_{j=1,\dots,2n-1 \atop k=1,\dots,n-1} \quad
[j(\lambda_k^j - \lambda_k^{-j})]_{j=1,\dots,2n-1 \atop k=1,\dots,n}
\Big ] \nonumber \\
&& \qquad = \prod_{j=1}^n (\lambda_j - 1/\lambda_j)
\prod_{1 \le j < k \le n} (\lambda_k - \lambda_j)^2
(1/\lambda_k - 1/\lambda_j)^2 (\lambda_j - 1/\lambda_k)^2
(1/\lambda_j - \lambda_k)^2.
\end{eqnarray}
\end{prop}

\noindent
Proof. \quad We see that the determinant is a symmetric rational function
in $\lambda_1,\dots,\lambda_n$, and is antisymmetric under the mapping
$\lambda_i \mapsto 1/\lambda_i$ for any $i=1,\dots,n$. It must thus
be of the form
\begin{equation}\label{Q1}
\prod_{j=1}^n (\lambda_j - 1/\lambda_j) \,
q(\lambda_1,\dots,\lambda_n)
\end{equation}
where $q$ is symmetric and unchanged by the mapping $\lambda_i 
\mapsto 1/\lambda_i$. Noting too that the determinant vanishes when
$\lambda_1 = \lambda_2$, we see that $q$ must contain as a factor
\begin{equation}\label{Q2}
\prod_{1 \le j < k \le n} (\lambda_k - \lambda_j)^2(1/\lambda_k -
1/\lambda_j)^2 (\lambda_j - 1/\lambda_k)^2 (1/\lambda_j - \lambda_k)^2.
\end{equation}
The highest order term (in degree) of (\ref{Q2}) multiplied by
$\prod_{j=1}^n(\lambda_j - 1/\lambda_j)$ is
$\prod_{j=1}^n\lambda_j \prod_{1 \le j < k \le n}
(\lambda_k - \lambda_j)^4$. On the other hand the highest order term
in degree in the determinant is
$$
\det \Big [ [\lambda_k^j - \lambda_n^j]_{j=1,\dots,2n-1 \atop
k = 1,\dots, n-1}
[j \lambda_k^{j}]_{j=1,\dots,2n-1 \atop k=1,\dots,n} \Big ] .
$$
According to (\ref{c}) this evaluates to the same expression, so in fact
$q$ must be exactly equal to (\ref{Q2}).
\hfill $\square$

\subsection{The Jacobians}
Using a similar strategy to that used to derive the Jacobian (\ref{qb}) in
the proof of Theorem \ref{th1}, the results of the previous subsection
together with the method of wedge products allows the two Jacobians
evaluated by indirect means in \cite{KN04} to be derived directly.

\begin{thm}
Consider unitary Hessenberg matrices with entries specified by (\ref{Hr})
and surrounding text. The Jacobian for the change of variables from
$\{\alpha_j\}_{j=0,\dots,n-1}$ to $\{\theta_i\}_{i=1,\dots,n}$,
$\{q_i\}_{i=1,\dots,n-1}$ is equal to
\begin{equation}\label{a.1}
{\prod_{i=0}^{n-2}(1 - |\alpha_i|^2) \over q_n \prod_{i=1}^n q_i}.
\end{equation}

Consider $2n \times 2n$ real orthogonal Hessenberg matrices as specified
above Proposition \ref{I3}. The Jacobian for the change of variables from
$\{\alpha_j\}_{j=0,\dots,2n-2}$ to $\{\theta_i\}_{i=1,\dots,n}$,
$\{q_i\}_{i=1,\dots,n-1}$ is equal to
\begin{equation}\label{a.2}
{2^{n-1} \over q_n \prod_{i=1}^n q_i}
{\prod_{l=0}^{2n-2}(1 - |\alpha_l|^2) \over
\prod_{k=0}^{2n-2}(1 - \alpha_k)^{1/2}(1 + (-1)^k \alpha_k)^{1/2} }.
\end{equation}
\end{thm}

\noindent
Proof. \quad In relation to (\ref{a.1}) we begin by equating successive powers
of $\lambda$ on both sides of (\ref{aH}). Recalling the explicit form of $H$
given by (\ref{Hr}) and surrounding text this gives
\begin{eqnarray}\label{gipz}
1 & = & \sum_{j=1}^n q_j^2 \nonumber \\
\bar{\alpha}_0 & = & \sum_{j=1}^n q_j^2 \lambda_j \nonumber \\
* + \bar{\alpha}_1 \rho_0^2 & = & \sum_{j=1}^n q_j^2 \lambda_j^2 
 \nonumber \\
* + \bar{\alpha}_2 \rho_0^2 \rho_1^2 & = & \sum_{j=1}^n q_j^2 \lambda_j^3
 \nonumber \\
\vdots & & \vdots \nonumber \\
* + \bar{\alpha}_{n-1} \rho_0^2 \rho_1^2 \cdots \rho_{n-2}^2 & = &
\sum_{j=1}^n q_j^2 \lambda_j^n
\end{eqnarray}
where the $*$ denotes terms involving only variables already having
appeared on the l.h.s.~of the preceding equation.  Thus as in the 
corresponding equations (\ref{2.q})
in the tridiagonal case a triangular
structure results. We know that $\alpha_j$ $(j=0,\dots,n-2)$ has an
independent real and imaginary part, while $\alpha_{n-1}:=
e^{i \phi}$, $\lambda_j := e^{i \theta_j}$ $(j=1,\dots,n)$ have unit
modulus. Consequently the number of equations can be made equal to the
number of variables by firstly using the first equation to eliminate
$q_n^2$ in the subsequent equations, then appending to the list the complex
conjugate of all but the last of the remaining equations.
 
Let us take differentials of these $2n-1$ equations, then take wedge
products of both sides. Because of the triangular
structure, we obtain on the l.h.s. 
\begin{equation}\label{C.4}
\rho_0^2 \rho_1^2 \cdots \rho_{n-2}^2 \prod_{l=1}^{n-2} \rho_{n-l-2}^{4l}
d \vec{\alpha} \wedge d \phi,
\end{equation}
while this operation on the r.h.s.~yields
\begin{eqnarray}\label{C.5}
&&q_n^2 \prod_{j=1}^{n-1} q_j^3 \left |
\det \left [ \begin{array}{cc}
\left [ \begin{array}{c}\lambda_k^j - \lambda_n^j  \\
\lambda_k^{-j} - \lambda_n^{-j} \end{array}
\right ]_{j,k=1,\dots,n-1}
&
\left [ \begin{array}{c}j \lambda_k^j \\
-j \lambda_k^{-j}
\end{array} \right ]_{j=1,\dots,n-1 \atop k=1,\dots,n}
\\
{}[\lambda_k^n - \lambda_n^n]_{k=1,\dots,n-1} &
[n \lambda_k^n]_{k=1,\dots,n}
\end{array} \right ] \right |
d\vec{\theta} \wedge d \vec{q} \nonumber \\
&& \qquad = q_n^2 \prod_{j=1}^{n-1} q_j^3 \prod_{1 \le j < k \le n}
| \lambda_k - \lambda_j|^4 d\vec{\theta} \wedge d \vec{q}
\end{eqnarray}
where the equality follows upon using the determinant evaluation
(\ref{TN}).

Analogous to (\ref{2.14'}), by definition the Jacobian $J$ satisfies
$$
d \vec{\alpha} \wedge d \phi = J d \vec{\theta} \wedge d \vec{q}.
$$
Equating (\ref{C.4}) and (\ref{C.5}) and making use of (\ref{3H}) gives
(\ref{a.1}).

Consider next the derivation of (\ref{a.2}). Proceeding as in the
derivation of (\ref{gipz}), expanding (\ref{3.22'}) in powers of $\lambda$,
we obtain 
\begin{eqnarray*}
1 & = & \sum_{j=1}^n q_j^2 \\
\alpha_0 & = & {1 \over 2} \sum_{j=1}^n q_j^2 (\lambda_j + \bar{\lambda}_j)
\\
* + \alpha_1 \rho_0^2 & = & {1 \over 2} \sum_{j=1}^n q_j^2
(\lambda_j^2 + \bar{\lambda}_j^2) \\
\vdots & & \vdots \\
* + \alpha_{2n-2} \rho_0^2 \cdots \rho_{2n-3}^2 & = & {1 \over 2}
\sum_{j=1}^n q_j^2 (\lambda_j^{2n-1} + \bar{\lambda}_j^{2n-1}).
\end{eqnarray*}
The l.h.s.~again exhibits a triangular structure, and furthermore all
quantities on the l.h.s.~are real. Taking the differentials of both sides,
and forming the wedge product of the l.h.s.'s of all but the first equation
gives
\begin{equation}\label{d.x}
\prod_{l=0}^{2n-3} \rho_l^{2(2n-2-l)} d \vec{\alpha}.
\end{equation}
On the r.h.s., after substituting for $q_n d q_n$ using the differential of
the first equation, this same procedure gives
\begin{equation}\label{d.x1}
2^{-n} q_n^2 \prod_{j=1}^{n-1} q_j^3 
\left | 
\det \Big [ [ \lambda_k^j + \lambda_k^{-j} - (\lambda_k^j + \lambda_k^{-j})
]_{j=1,\dots,2n-1 \atop k=1,\dots,n-1} \quad
[j(\lambda_k^j - \lambda_k^{-j})]_{j=1,\dots,2n-1 \atop k=1,\dots,n}
\Big ] \right | d \vec{\theta} \wedge d \vec{q}.
\end{equation}
Here the Jacobian $J$ satisfies
$$
d \vec{\alpha} = J d \vec{\theta} \wedge d \vec{q},
$$
so (\ref{d.x}) and (\ref{d.x1}) (with the determinant evaluated according to
(\ref{Q0})) together give a formula for $J$ in terms of $\{q_i\}$,
$\{ \alpha_i \}$ and $\{\lambda_i \}$. The latter set of variables can be
eliminated by making use of Proposition \ref{I3}, and
(\ref{a.2}) results.
\hfill $\square$

\section{A multiplicative rank 1 perturbation of unitary matrices}
\setcounter{equation}{0}
\subsection{Circular analogue of the Dixon-Anderson density}\label{s4.1}
Let $\vec{e}_1$ denote the $n \times 1$ unit vector $(1,0,\dots,0)^T$.
Let $t$ be a complex number with $|t|=1$. Then the matrix
$$
I_n - (1-t) \vec{e}_1 \vec{e}_1^T
$$
is a unitary matrix differing from the identity only in the top left entry
which is $t$. Our interest in this section is the
eigenvalue distribution of
\begin{equation}\label{IU}
\tilde{U} := (I_n - (1-t) \vec{e}_1 \vec{e}_1^T) U,
\end{equation}
for $U$ a given unitary matrix. Such multiplicative rank 1 perturbations
are discussed for example in \cite{AGR91}. 
The term multiplicative perturbation is used because
$\tilde{U}$ is obtained
from $U$ by multiplication of the first row by the unimodular complex
number $t$, while the term rank 1 is used because the
multiplicative perturbative factor
differs from the identity by a rank 1 matrix.
We will see that for $U$ a member of the circular
$\beta$-ensemble, a joint eigenvalue p.d.f.~generalizing (\ref{1.2})
results.

First a rational function having as its zeros the eigenvalues of
$\tilde{U}$ will be specified.

\begin{prop}\label{p6}
Let $U$ be an $n \times n$ unitary matrix with distinct eigenvalues
$e^{i \theta_1}, \dots, e^{i \theta_n}$, and denote the corresponding
matrix of eigenvectors by $V = [v_{jk}]_{j,k=1,\dots,n}$. The
eigenvalues of $\tilde{U}$ as specified by (\ref{IU}) occur at the
zeros of the rational function
\begin{equation}\label{et}
C_n(\lambda) = 1 + (t-1) \sum_{j=1}^n {e^{i \theta_j} |v_{1j}|^2 \over 
e^{i \theta_j} - \lambda}.
\end{equation}
\end{prop}

\noindent
Proof. \quad Noting from (\ref{IU}) that $\tilde{U}=U - (1-t) U'$, where
$U'$ is the $n \times n$ matrix in which the first row is equal to the
first row of $U$, and all other rows have all entries zero, we see
$$
V^{-1} \tilde{U} V = {\rm diag}[e^{i\theta_1},\dots,e^{i \theta_n}]
+ (t-1)[\bar{v}_{1j} v_{1k} e^{i\theta_k} ]_{j,k=1,\dots,n}.
$$
Thus $\tilde{U}$ has the same spectrum as a matrix which consists of a
rank 1 multiplicative perturbation of a diagonalized unitary matrix.
The characteristic polynomial of this matrix can be factorized as
$$
\prod_{l=1}^n(e^{i \theta_l} - \lambda)
\det \Big [ \delta_{j,k} + (t-1) \bar{v}_{1j} v_{1k} e^{i \theta_k}/
(e^{i \theta_j} - \lambda) \Big ]_{j,k=1,\dots,n},
$$
and the zeros must occur at the zeros of the determinant. Noting the
simple determinant evaluation
$$
\det [u_j \delta_{j,k} + 1]_{j,k=1,\dots,n} =
\prod_{l=1}^n u_l \Big ( 1 + \sum_{j=1}^n {1 \over u_j} \Big ),
$$
the sought result follows.
\hfill $\square$

\medskip
We remark that Proposition \ref{p6} can be extended to the case that each
eigenvalue $e^{i \theta_j}$ has multiplicity $m_j$. Thus with
$v_{1j}^{(s)}$ denoting the first component of the $s$th independent eigenvector
corresponding to $e^{i \theta_j}$, we replace $|v_{1j}|^2$ in
(\ref{et}) by $\sum_{s=1}^{m_j} |v_{1j}^{(s)}|^2$.

Let us suppose now that the matrix $U$ is a unitary upper triangular
Hessenberg matrix parametrized as specified by (\ref{Hr}) and 
surrounding text. One of the main results of \cite{KN04} is that the
parameters $\{\alpha_j\}_{j=0,\dots,n-1}$ can be chosen from particular
probability distributions so that the eigenvalue p.d.f.~of $U$ is given
by (\ref{1.2}). 
The probability distributions in question are 
parametrized by a real number $\nu \ge 1$ and denoted by $\Theta_\nu$.
For $\nu > 1$, the support of $\Theta_\nu$ is the open unit disk 
$|z|<1$
in 
the complex plane, and the distribution is specified by the
p.d.f.
$$
{\nu -1 \over 2 \pi} (1 - |z|^2)^{(\nu - 3)/2}.
$$
For $\nu = 1$, the support is the unit circle $|z|=1$, and
$\Theta_1$ denotes the uniform distribution. Proposition 4.2 of
\cite{KN04} tells us that if
\begin{equation}\label{he}
\alpha_{n-j-1} \sim \Theta_{\beta j + 1} \qquad (j=0,\dots,n-1)
\end{equation}
then the corresponding eigenvalue p.d.f.~is given by (\ref{1.2}). 
Furthermore, it tells us that the modulus squared of the first
component of the eigenvectors $|v_{1j}|^2 := \mu_j$ have the distribution
with measure
\begin{equation}\label{4.24'}
{1 \over C_{\beta N} } \prod_{i=1}^n \mu_i^{\beta/2 - 1} d \vec{\mu},
\end{equation}
where
$$
C_{\beta N} = {\Gamma^N(\beta/2) \over \Gamma(\beta N/2) }, \qquad
d \vec{\mu} := d\mu_1 \dots d\mu_{n-1}.
$$
This is an example of the Dirichlet distribution.

The latter fact motivates the study of the zeros of (\ref{et}) with
the $|v_{1j}|^2$ distributed according to the Dirichlet distribution.
We will find that a conditional p.d.f.~relating to (\ref{1.2}) results
provided the distribution of $t$ is appropriately chosen. First, some
preliminary results must be established.

\begin{lemma}\label{le1}
Suppose in (\ref{et}) that
\begin{equation}\label{tha1}
0 < \theta_1 < \theta_2 < \cdots < \theta_n \le 2 \pi.
\end{equation}
The function $C(\lambda)$ has exactly $n$ zeros occurring at
$\lambda = e^{i \psi_1}, \dots, e^{i \psi_n}$, where
\begin{equation}\label{tha}
\theta_{i-1} < \psi_i < \theta_i \qquad (i=1,\dots,n, \quad
\theta_0 := \theta_n {\rm mod} \, 2 \pi).
\end{equation}
Furthermore, with $\lambda_j := e^{i \theta_j}$,
$\tilde{\lambda}_j := e^{i \psi_j}$, we have
\begin{eqnarray}
- (t-1) \lambda_j q_j & = & {\prod_{l=1}^n(\lambda_j - \tilde{\lambda}_l)
\over \prod_{l=1,l\ne j}^n(\lambda_j - \lambda_l) } \qquad
(j=1,\dots,n) \label{r1} \\
\prod_{l=1}^n \tilde{\lambda}_l & = & t \prod_{l=1}^n \lambda_l. \label{r0}
\end{eqnarray}
\end{lemma}

\noindent
Proof. \quad The fact that there are exactly $n$ zeros of unit modulus
follows from the relationship of $C_n(\lambda)$ to the characteristic
polynomial of a unitary matrix. The interlacing property is well known
\cite{AGR91}. It can be seen graphically by writing (\ref{et}) in the form
$$
C_n(\lambda) = {(t-1) \over 2i} \Big (
\cot {\phi \over 2} - \sum_{i=1}^n q_i \cot \Big ( {\psi - \theta_i \over
2} \Big ) \Big ),
$$
where we have set $t:= e^{i \phi}$, $\lambda := e^{i \psi}$.

With the zeros so identified, and the poles as evident from (\ref{et}),
it follows that
\begin{equation}\label{et1}
C_n(\lambda) = {\prod_{j=1}^n(\lambda - \tilde{\lambda}_j) \over
\prod_{j=1}^n (\lambda - \lambda_j) },
\end{equation}
where use has also been made of the property $C_n(\lambda) \to 1$ as
$\lambda \to \infty$. Comparing residues in (\ref{et}) and (\ref{et1})
gives (\ref{r1}), while setting $\lambda = 0$ gives (\ref{r0}).
\hfill $\square$

\medskip
The Jacobians for some change of variables are also required.

\begin{lemma}\label{le2}
Let $J$ be the Jacobian for the change of variables
$\{q_j\}_{j=1,\dots,n-1} \cup \{\phi \}$ to $\{\psi_j \}_{j=1,\dots,n}$.
We have
\begin{equation}\label{4.6}
J  = 
 =  |1 - t|^{-(n-1)}
\prod_{1 \le j < k \le n} \Big |
{\tilde{\lambda}_k - \tilde{\lambda}_j \over
\lambda_k - \lambda_j} \Big |
\end{equation}
\end{lemma}

\noindent
Proof. \quad By definition $J$ is positive and satisfies
\begin{equation}\label{4.14'}
d \vec{q} \wedge dt = J d \vec{\psi}.
\end{equation}
Now
\begin{eqnarray*}
d \vec{q} \wedge d \phi & = & (t-1)^{-(n-1)} d(t-1) \vec{q} \wedge
d \phi \nonumber \\
& = & (t-1)^{-(n-1)} \det \bigg [ \Big [ {\partial (t-1) q_l \over \partial
\tilde{\lambda}_j } \Big ]_{j=1,\dots,n \atop l=1,\dots,n-1}
\Big [ {\partial t \over \partial \tilde{\lambda}_j }
\Big ]_{j=1,\dots,n} \bigg ]
d \vec{\tilde{\lambda}}_j.
\end{eqnarray*}
But according to (\ref{r1}) and (\ref{r0})
$$
{\partial (t-1) q_l \over \partial \tilde{\lambda}_j }
= {(t-1) q_l \over \lambda_l - \tilde{\lambda}_j }, \qquad
{\partial t \over \partial \tilde{\lambda}_j } =
{t \over \tilde{\lambda}_j },
$$
and so with $\lambda_n = 0$ (temporarily as a notational convenience)
\begin{equation}\label{lh}
d \vec{q} \wedge d \phi =
t \prod_{l=1}^{n-1} (-q_l)
\det \Big [ {1 \over \tilde{\lambda}_j - \lambda_l} \Big ]_{j,l=1,\dots,n}
d \vec{\tilde{\lambda}}.
\end{equation}
Since $J$ is positive and satisfies (\ref{4.14'}), it must be equal to the
modulus of the terms multiplying $d \vec{\tilde{\lambda}}$ in this expression.
Evaluating the determinant as a Cauchy double alternant, and evaluating
$\prod_{l=1}^{n-1} q_l $ using (\ref{r1}) gives the stated result.
\hfill $\square$

\medskip
The results of the above two lemmas allow a change of variables to
be made from $\{q_j\}_{j=1,\dots,n-1} \cup \{t\}$ to
$\{\lambda_l \}_{l=1,\dots,n}$.

\begin{thm}\label{tm5}
With $|v_{1j}|^2 = q_j$ $(j=1,\dots,n)$ in (\ref{et}), let $\{q_j\}$ have
the Dirichlet distribution with measure
\begin{equation}\label{4.4}
{\Gamma((n-1)d + d_0) \over (\Gamma(d))^{n-1} \Gamma(d_0) }
\Big ( \prod_{j=1}^{n-1} q_j^{d-1} \Big ) q_n^{d_0 - 1} d \vec{q}.
\end{equation}
Further, let the parameter $t$ in (\ref{et}) be determined by the p.d.f.~with
measure
\begin{equation}\label{4.5}
{\Gamma^2({1 \over 2} (d_0 + (n-1)d + 1)) \over
2 \pi \Gamma((n-1)d + d_0)}
|1 - t|^{d_0 + (n-1)d - 1} d \phi.
\end{equation}
The conditional p.d.f.~of $\{\tilde{\lambda}_j = e^{i \psi_j}
\}_{j=1,\dots,n}$, given $\{\lambda_j = e^{i \theta_j} \}_{j=1,\dots,n}$,
is equal to
\begin{equation}\label{3.Ap}
A {\prod_{l=1}^n|e^{i \theta_n} - e^{i \psi_l}|^{d_0-1} \over
\prod_{l=1}^{n-1}|e^{i \theta_n} - e^{i \theta_l}|^{d_0 + d - 1}}
{\prod_{j=1}^{n-1} \prod_{l=1}^n |e^{i \theta_j} - e^{i \psi_l}
|^{d - 1}
\over \prod_{1 \le j < k \le n-1} | e^{i \theta_k} - e^{i \theta_j}
|^{2 d - 1}}
\prod_{1 \le j < k \le n} |e^{i \psi_k} - e^{i \psi_j}|,
\end{equation}
$$
A := {\Gamma((n-1)d + d_0) \over (\Gamma(d))^{n-1}
\Gamma(d_0)}
{\Gamma^2({1 \over 2} (d_0 + (n-1)d + 1)) \over
2 \pi \Gamma((n-1)d + d_0)}.
$$
\end{thm}

\noindent
Proof. \quad Our proof, which at a technical level proceeds by
making use of the results of Lemmas \ref{le1} and \ref{le2}, is based
on a strategy adopted for an analogous problem with real roots
by Anderson \cite{An91}, and many years before by Dixon \cite{Di05}.

We must change variables in the product of (\ref{4.4}), (\ref{4.5}),
and the Jacobian $J$. We know the latter is given by  (\ref{4.6}).
To change variables in (\ref{4.4}) we use  (\ref{r0}), which gives
\begin{eqnarray}\label{28}
&& (\prod_{j=1}^{n-1} q_j^{d-1} ) q_n^{d_0 - 1} =
{1 \over |1 - t|^{d_0 +(n-1) d - n}}
{ \prod_{j=1}^{n-1} \prod_{l=1}^n |\lambda_j - \tilde{\lambda}_l|^{d-1}
\over \prod_{1 \le j < l \le n-1} |\lambda_l - \lambda_j|^{2(d-1)}}
{\prod_{l=1}^n|\lambda_n - \tilde{\lambda}_l|^{d_0-1} \over
\prod_{l=1}^{n-1} |\lambda_n - \lambda_l|^{d_0+d - 2}}.
\end{eqnarray}
Multiplying (\ref{4.5}), (\ref{4.6}) and (\ref{28}) we see that the 
dependence on $t$ cancels, and the expression (\ref{3.Ap}) results.
\hfill $\square$

\subsection{Properties of the corresponding joint density}\label{s4.2}
Let us suppose $\{\theta_j\}_{j=1,\dots,n}$, assumed ordered as in
(\ref{tha}) and with $\theta_n$ fixed, are distributed according to the p.d.f.
$$
{(n-1)! \over (2 \pi)^{n-1} M_{n-1}((a_1+d_0 + d -1)/2,
(a_1+d_0 + d - 1)/2,d) }
\prod_{l=1}^{n-1} |e^{i \theta_n} - e^{i \theta_l} |^{a_1 + d_0 + d - 1}
\prod_{1 \le j < k \le n - 1} |e^{i \theta_k} - e^{i \theta_j} |^{2d},
$$
\begin{equation}\label{4.16'}
M_N(a,b,\lambda) =
\prod_{j=0}^{N - 1} { \Gamma (\lambda j+a+b+1)
\Gamma(\lambda (j+1)+1) \over
  \Gamma (\lambda j+a+1)\Gamma (\lambda j+b+1) \Gamma (1 + \lambda)},
\end{equation}
(for a discussion of this p.d.f.~see \cite{Fo02}). Multiplying this with
(\ref{3.Ap}) gives the joint p.d.f.
\begin{eqnarray}\label{cs}
&&C^{(n,n-1)}(\psi,\theta) := {A (n-1)! \over
(2 \pi )^{n-1} M_{n-1}((a+a_1+d)/2,(a+a_1+d)/2,d)}
\prod_{l=1}^n|e^{i \theta_n} - e^{i \psi_l}|^{a}
\nonumber \\
&& \qquad \times
\prod_{1 \le j < k \le n} |e^{i \psi_k} - e^{i \psi_j}|
\prod_{l=1}^{n-1}|e^{i \theta_n} - e^{i \theta_l}|^{a_1}
\prod_{1 \le j < k \le n-1} | e^{i \theta_k} - e^{i \theta_j}|
\prod_{j=1}^{n-1} \prod_{l=1}^n |e^{i \theta_j} - e^{i \psi_l}
|^{d - 1}
\end{eqnarray}
where $d_0 - 1 =: a$. 

The case of (\ref{cs}) relevant to the circular $\beta$-ensemble of
Killip and Nenciu is $a=d-1$, $a_1 = 1$ and $d=\beta/2$. Then (\ref{cs}) is
symmetric in $\{\theta_l \}_{l=1,\dots,n}$ and 
in $\{\psi_l \}_{l=1,\dots,n}$, and $\theta_n$ may again
be considered as variable ((\ref{cs}) should then be multiplied by
$n/2\pi$ to get the correct normalization). 
It corresponds to the joint eigenvalue 
p.d.f.~of a unitary Hessenberg matrix with parameters distributed according
to (\ref{he}), and thus with eigenvalue p.d.f.~(\ref{1.2}), and the same
unitary Hessenberg matrix perturbed by multiplication of the first row
by $t$. The factor $t$ is to be distributed according to (\ref{4.5}) with
$d_0 = d = 
\beta/2$. We know that the p.d.f.~for $\{\theta_l\}_{l=1,\dots,n}$ can
be sampled by computing the zeros of $\chi_n(\lambda)$ as calculated
from (\ref{1.3}) with $\{\alpha_j\}_{j=0,\dots,n-1}$ chosen as specified
by (\ref{he}). To sample from $\{\psi_l\}_{l=1,\dots,n-1}$ in the joint
p.d.f., with the same $\{\alpha_j\}_{j=0,\dots,n-1}$ we again compute
$\chi_n(\lambda)$ from the recurrences (\ref{1.3}), but now with
$\chi_0(\lambda) = \tilde{\chi}_0 = t$. 

Next, let us turn our attention to integration formulas associated with
(\ref{cs}). Since (\ref{3.Ap}) is a conditional p.d.f.~we must have
\begin{eqnarray}\label{129}
&&
\int_R d \psi_1 \cdots d \psi_n \, C^{(n,n-1)}(\psi, \theta) =
 {(n-1)! \over (2 \pi )^{n-1}
M_{n-1}((a+a_1+d)/2,(a+a_1+d)/2,d)}
\\
&& \qquad \times
\prod_{l=1}^{n-1}|e^{i \theta_n} - e^{i \theta_l}|^{a+a_1+d}
\prod_{1 \le j < k \le n-1} | e^{i \theta_k} - e^{i \theta_j}|^{2 d}
\end{eqnarray}
where $R$ denotes the region specified by the inequalities (\ref{tha1}).
Special cases of (\ref{129}) are two classical inter-relations between
circular ensembles \cite{Dy62a,MD63} (for an extensive study of
such formulas in random matrix theory see \cite{FR01}, and for application
of the Dixon-Anderson density to the cases with real eigenvalues
see \cite{FR02b}).

\begin{prop}\label{CC}
Let COE${}_n$, CUE${}_n$, CSE${}_n$ --- the circular ensembles with
orthogonal, unitary and symplectic symmetry respectively --- refer to the
eigenvalue p.d.f.~(\ref{1.2}) with $\beta=1,2,4$ respectively. Let alt
refer to the operation of integrating over every second eigenvalue. 
Let ${\rm COE}{}_n \cup {\rm COE}{}_n$ denote the p.d.f.~of $2n$
eigenvalues which results from superimposing two independent sequences of
$n$ eigenvalues each with a COE${}_n$ distribution. 
One has
\begin{eqnarray}
{\rm alt}( {\rm COE}{}_{n} \cup {\rm COE}{}_n) & = & {\rm CUE} {}_n \label{t1}
\\{\rm alt}( {\rm COE}{}_{2n}) & = & {\rm CSE} {}_n \label{t2} 
\end{eqnarray}
\end{prop}

\noindent
Proof. \quad For the first identity we require the fact that \cite{Gu62}
\begin{equation}\label{s2}
{\rm COE}_n \cup {\rm COE}_n \propto
\prod_{1 \le j < k \le n}
|e^{i \theta_{2k}} - e^{i \theta_{2j}}|
|e^{i \theta_{2k-1}} - e^{i \theta_{2j-1}}|.
\end{equation}
We then
see that (\ref{129}) with $a_1 = d = 1$ is equivalent to the first
identity. The second identity is immediately seen to correspond to
(\ref{129}) with $a_1=2$, $d=2$.
\hfill $\square$

\subsection{Matrix theoretic derivation of the COE, CUE, CSE
inter-relations}
The inter-relations (\ref{t1}), (\ref{t2}) were originally proved
by establishing the same integration formulas as those noted in the
proof of Proposition  \ref{CC}. In this subsection it will be shown
how random matrices can be constructed in such a way that both
(\ref{t1}) and (\ref{t2}) are immediate.

Let us consider first (\ref{t1}). This requires a different random
matrix realization of the joint p.d.f.~(\ref{cs}) in the case
$d=a_1=1$, $a=0$ to that given in the paragraph below (\ref{cs}).
The theory underlying the construction is the following.

\begin{prop}\label{ps1}
Let $M_1$ be a $2n \times 2n$ unitary matrix with real elements, which has
a doubly degenerate spectrum with the independent eigenvalues distributed
as CUE${}_n$. Let the matrix of eigenvectors be $V=[v_{ij}]_{i,j=1,\dots,
2n}$, and suppose the joint distribution of $\mu_j := (v_{1\, 2j-1})^2 +
(v_{1 \, 2j})^2$ $(j=1,\dots,n)$ is equal to the Dirichlet distribution
(\ref{4.24'}) with $\beta = 2$. Form the matrix $M_1'$ by multiplying the
first row of $M_1$ by the complex number $t$, $|t|=1$, where $t$ has
distribution (\ref{4.5}) with $d_0=d=1$. Then the perturbed matrix
$M_1'$ has for its eigenvalue p.d.f.~(\ref{cs}) with $a=0$, $a_1=d=1$.
\end{prop}

\noindent
Proof. \quad Let $\{e^{i \theta_l}\}_{l=1,\dots,n}$ denote the
independent eigenvalues of $M_1$. Proceeding as in the derivation of
(\ref{et}) shows that the characteristic polynomial of the
perturbed matrix is equal to
\begin{equation}\label{4.23'}
\prod_{j=1}^n(e^{i \theta_l} - \lambda)^2 \Big (
1 + (t-1) \sum_{j=1}^n {e^{i \theta_j} \mu_j \over
e^{i \theta_j} - \lambda} \Big ).
\end{equation}
Thus $M_1'$ has $n$ eigenvalues at $\{e^{i \theta_l}\}_{l=1,\dots,n}$, and
$n$ eigenvalues given by the zeros of the rational function factor in
(\ref{4.23'}). We are given that the former have p.d.f.~CUE${}_n$, while 
Theorem \ref{tm5} tells us that the latter have conditional
p.d.f.~(\ref{3.Ap}) with $d_0=d=1$. Multiplying these together gives
the stated joint distribution.
\hfill $\square$

\medskip
To realize the matrix $M_1$, we begin with an element of $U(n)$ chosen
according to the Haar measure. This gives the eigenvalue p.d.f.~CUE${}_n$,
with the eigenvectors such that the $\mu_j := |v_{1j}|^2$ $(j=1,\dots,n)$
have joint distribution (\ref{4.24'}) with $\beta = 2$. To obtain a
doubly degenerate spectrum, each element $x+iy$ is replaced by its
$2 \times 2$ real matrix representation
$$
\left [ \begin{array}{cc} x & y \\ -y & x \end{array} \right ],
$$
so forming a $2n \times 2n$ matrix with real entries.
Since the corresponding perturbed matrix $M_1'$ retains all distinct
eigenvalues of $M_1$,
\begin{equation}\label{pM}
{\rm unpert} (M_1') = {\rm CUE}_n,
\end{equation}
where the l.h.s.~denotes the eigenvalue p.d.f.~of $M_1'$ integrated over the
perturbed eigenvalues. On the other hand Proposition \ref{ps1} together
with (\ref{s2}) tell us that with reference to the eigenvalue
p.d.f., $M_1' = {\rm COE}_n \cup {\rm COE}_n$. Thus we have a matrix
theoretic
understanding of (\ref{t1}) in the sense that its validity is a 
consequence of spectral properties of $M_1'$ which avoid the need for
explicit integration of the eigenvalue p.d.f.

We seek a similar understanding of (\ref{t2}). For this we must identify an
ensemble of random matrices with a doubly degenerate spectrum, and
their perturbations, which give rise to (\ref{cs}) in the
case $d=2$, $a=a_1=1$. In fact the very definition of the circular
symplectic ensemble involves matrices with a doubly degenerate spectrum
(see e.g.~\cite{Fo02}). Thus, if for any $2n \times 2n$ matrix $X$
we set
$$
X^D := Z_{2n} X^T Z_{2n}, \qquad {\rm where} \quad
Z_{2n} := I_n \otimes \left [ \begin{array}{cc}0 & -1 \\ 1 & 0
\end{array} \right ],
$$
and select $U \in U(2n)$ with Haar measure, then matrices of the form
$U^DU$ make up the circular symplectic ensemble. Such matrices have a
doubly degenerate spectrum, and the $n$ independent eigenvalues are
distributed according to CSE${}_n$. Furthermore, with the matrix of
eigenvectors denoted $V= [v_{ij}]_{i,j=1,\dots,2n}$, one has that the
$\mu_j : = |v_{1 \, 2j-1}|^2 + |v_{1 \, 2j}|^2$ $(j=1,\dots,n)$ are
distributed according to the Dirichlet distribution (\ref{4.24'}) with
$\beta = 4$. Consideration of these facts, together with reasoning
analogous to that used in the proof of Proposition \ref{ps1}, gives the
sought realization.

\begin{prop}\label{p8}
Let $M_2$ be a member of the circular symplectic ensemble as specified
above. Form the matrix $M_2'$ by multiplying the first row of $M_2$
by the complex number $t$, $|t|=1$, where $t$ has distribution
(\ref{4.5}) with $d_0=d=2$. The joint eigenvalue p.d.f.~of $M_2'$ is
then given by (\ref{cs}) with $a=a_1=1$, $d=2$.
\end{prop}

Analogous to (\ref{pM}) it is immediate that
$$
{\rm unpert}(M_2') = {\rm CSE}{}_n.
$$
Because Proposition \ref{p8} tells us that $M_2'$ has a joint distribution
formally equivalent to COE${}_{2n}$, (\ref{t2}) is reclaimed as a matrix
theoretic identity.

\section{Cayley transformation}
\subsection{Cauchy analogue of the Dixon-Anderson density}
\setcounter{equation}{0}
In general a unitary matrix $U$ is transformed to an Hermitian matrix $H$ by
the Cayley transformation
\begin{equation}\label{HU}
H = i {1_N - U \over 1_N + U}.
\end{equation}
At the level of the eigenvalues, the change of variables (\ref{HU}) in the workings
of Sections \ref{s4.1} and \ref{s4.2}  leads to a joint p.d.f.~on interlacing
variables on the real line, relating to the so called (generalized) Cauchy
ensemble \cite{WF01,BO01}. Properties of this allow a random three term recurrence
to be derived for the (projected)
characteristic polynomial associated with the p.d.f.~(\ref{Mn}).

First we apply the change of variables implied by (\ref{HU}) to (\ref{et}) with the
l.h.s.~written as (\ref{et1}).

\begin{prop}
Consider the rational function (\ref{et}) with the lower terminal of summation
extended to $j=0$. Substitute for the l.h.s.~(\ref{et1}) with the
lower terminals in the products extended to $j=0$. Under the change of variables
\begin{eqnarray}\label{cc1}
&& \tilde{\lambda}_j = {x_j - i \over x_j + i}, \qquad
\lambda_j = {y_j - i \over y_j + i} \quad (j \ne 0) \nonumber \\
&& \lambda = {z - i \over z + i}, \qquad t = {c - i \over c + i}
\end{eqnarray}
and with $\lambda_0 = 1$ we obtain
\begin{equation}\label{cc2}
{\prod_{j=0}^n (z - x_j) \over (z^2+1) \prod_{j=1}^n (z - y_j) } =
{z - c \over q_0 (z^2 + 1) } - \sum_{j=1}^n {(q_j/q_0) \over z - y_j},
\end{equation}
where
\begin{equation}\label{cc2'}
x_0 > y_1 > x_1 > \cdots > y_n > x_{n+1}.
\end{equation}
\end{prop}

\noindent
Proof. \quad This follows from direct substitution, together with the formula
\begin{equation}\label{cc3}
{q_0 \over c + i} = {\prod_{l=1}^n (y_l + i) \over \prod_{l=0}^n (x_l + i) },
\end{equation}
which is a consequence of (\ref{r1}), with lower product terminals extended to $l=0$,
in the case $j=0$.
\hfill $\square$

\medskip
\begin{thm}
Consider the rational function (\ref{cc2}). Let $\{q_j\}_{j=0,\dots,n-1}$ have
the Dirichlet distribution with measure
\begin{equation}\label{w1}
{\Gamma(\sum_{j=0}^n d_j) \over \prod_{j=0}^n \Gamma(d_j)}
\prod_{j=0}^n q_j^{d_j - 1} \, d \vec{q}.
\end{equation}
Let $c$ have the generalized Cauchy distribution with measure
\begin{equation}\label{w2}
{\Gamma(\gamma) \Gamma(\bar{\gamma}) \over
\pi 2^{2(1 - {\rm Re} \,\gamma)} \Gamma(2 {\rm Re} \, \gamma - 1) }
(1 + i c)^{-\gamma}(1 - ic)^{-\bar{\gamma}} \, dc
\end{equation}
where
\begin{equation}\label{w3}
\sum_{i=0}^n d_i + 1 = 2 {\rm Re} \, \gamma.
\end{equation}
We have that the conditional p.d.f.~of $\{x_j\}_{j=0,\dots,n}$ given
$\{y_j\}_{j=1,\dots,n}$ is equal to
\begin{eqnarray}\label{w4}
&& \tilde{A} \prod_{j=0}^n (1 + i x_j)^{-\gamma} (1 - i x_j)^{-\bar{\gamma}}
\prod_{j=1}^n(1 + i y_j)^{\gamma - d_j} (1 - i y_j)^{\bar{\gamma} - d_j}
\nonumber \\
&& \qquad \times 
\prod_{j=1}^n \prod_{l=0}^n |y_j - x_l|^{d_j - 1}
\prod_{1 \le j < k \le n} |y_j - y_k|^{1 - d_j - d_k}
\prod_{0 \le j < k \le n} |x_j - x_k|
\end{eqnarray}
where
\begin{equation}\label{At}
\tilde{A} = {\Gamma(\gamma) \Gamma(\bar{\gamma}) \over
\pi 2^{2(1 - {\rm Re} \,\gamma}) }
{1 \over \Gamma(2 {\rm Re} \, \gamma - 1 - \sum_{i=1}^n d_i)
\prod_{j=1}^n \Gamma(d_j) }.
\end{equation}
\end{thm}

\noindent
Proof. \quad The task is to change variables in the wedge product of (\ref{w1})
and (\ref{w2}) to $\{x_j\}_{j=0,\dots,n}$. We have
\begin{eqnarray}\label{3.31}
(1 + i c)^{-\gamma} (1 - i c)^{-\bar{\gamma}}  & = &
|(1 - i c)^{-\bar{\gamma}} |^2 \nonumber \\
& = &
q_0^{-2 {\rm Re} \, \gamma}
\bigg | \bigg ( {\prod_{l=1}^n (1 - i y_l) \over
\prod_{l=0}^n (1 - i x_l) } \bigg )^{\bar{\gamma}} \bigg |^2 \nonumber \\
& = &
q_0^{-2 {\rm Re} \, \gamma}
{\prod_{l=1}^n (1 + i y_l)^\gamma (1 - i y_l)^{\bar{\gamma}} \over
\prod_{l=0}^n (1 + i x_l)^\gamma (1 - i x_l)^{\bar{\gamma}} }
\end{eqnarray}
where the second equality follows from (\ref{cc3}), and the final equality uses
the fact that since $x_l,y_l$ interlace according to (\ref{cc2'}),
$$
\log \bigg ( {\prod_{j=0}^n (1 + i x_j) \over \prod_{j=1}^n(1 + i y_j) }
\bigg ) = \sum_{j=0}^n \log (1 + i x_j) - \sum_{j=1}^n \log (1 + i y_j).
$$
Also, for $j=1,\dots,n$
\begin{equation}\label{3.322}
q_j = {q_0 \over |1 + i y_j|^2}
{\prod_{l=0}^n | y_j - x_l| \over \prod_{l=1,l \ne j}^n
|y_j - y_l| }.
\end{equation}

It remains to change variables in $d \vec{q} \wedge d c$. Since $x_j$ is
related to $\tilde{\lambda}_j$ and $c$ to $t=e^{i \phi}$
as given in (\ref{cc1}),
$$
d \vec{q} \wedge d c = {J \over 2}
|1 + i c|^2 \prod_{j=0}^n {2 \over (1 + i x_j) (1 - i x_j)}
d \vec{x}
$$
where $J$ is the Jacobian (\ref{4.6}) (appropriately modified to account for the
lower terminal being 0). In terms of the change of variables (\ref{cc1}) the
latter reads
$$
J = 2^{-n} q_0^n {\prod_{0 \le j < k \le n} |x_j - x_k|
\over \prod_{1 \le j < k \le n} |y_j - y_k| }
$$
and thus we have
\begin{equation}\label{3.45}
d \vec{q} \wedge dc =
q_0^{n+2} {1 \over |1 + i x_j|^2}
{\prod_{0 \le j < k \le n} |x_j - x_k|
\over \prod_{1 \le j < k \le n} |y_j - y_k| } \, d \vec{x}.
\end{equation}
Multiplying (\ref{3.31}), (\ref{3.322}) and (\ref{3.45}) gives the stated
result.
\hfill $\square$

\medskip
We remark that the conditional p.d.f.~(\ref{w4}) appears in \cite{Na03} as 
a generalization of a conditional p.d.f.~due to Dixon and
Anderson (see  (\ref{uk}) below). We remark too that the distribution
(\ref{w2}) in the case $\gamma$ real is the classical $t$-distribution.

\medskip
Integrating (\ref{w4}) over $\{x_j\}_{j=0,\dots,n}$ within the region
(\ref{cc2'}) we must get unity. Using this allows us to derive for the
multi-dimensional integral
\begin{equation}\label{IN}
I_n(\gamma;d) = {1 \over n!}
\int_{-\infty}^\infty dx_1 \cdots \int_{-\infty}^\infty dx_n \,
\prod_{l=1}^n (1 + i x_l)^{-\gamma} (1 - i x_l)^{-\bar{\gamma}}
\prod_{1 \le j < k \le n} |x_j - x_k|^{2 d}
\end{equation}
a recurrence analogous to that obtained by Anderson \cite{An91} for the
Selberg integral. Moreover, the intermediate workings will allow
us to deduce a random three term recurrence for the characteristic
polynomial associated with the p.d.f.~(\ref{Mn}).

\begin{cor}
We have
\begin{equation}\label{cf}
I_{n+1}(\gamma;d) = \pi 2^{2 - 2 {\rm Re} \, \gamma }
{\Gamma(2 {\rm Re} \, \gamma  - n d - 1) \Gamma((n+1)d) \over
\Gamma(d) | \Gamma(\gamma)|^2} I_n(\gamma-d;d).
\end{equation}
\end{cor}

\noindent
Proof. \quad Let us denote the region (\ref{cc2'}) by $R'$. As remarked, integrating
(\ref{w4}) over $\{x_j\}_{j=0,\dots,n}$ within $R'$ must give unity.
Setting
\begin{equation}\label{al}
d_1 = \cdots = d_n = d 
\end{equation}
this implies
\begin{eqnarray}\label{sgs}
&& \tilde{A}_d \int_{R'} dx_0 \cdots dx_n \,
\prod_{j=0}^n (1 + i x_j)^{-\gamma} (1 - i x_j)^{-\bar{\gamma}}
\prod_{j=1}^n \prod_{l=0}^n | y_j - x_l|^{d - 1}
\prod_{0 \le j < k \le n} |x_j - x_k| \nonumber \\
&& \qquad =
\prod_{j=1}^n(1 + i y_j)^{-\gamma + d}
(1 - iy_j)^{-\bar{\gamma} + d}
\prod_{1 \le j < k \le n} |y_j - y_k|^{2 d - 1}
\end{eqnarray}
where $\tilde{A}_d := \tilde{A} |_{d_1 = \cdots d_n = d}$. 
Thus
\begin{eqnarray}\label{66}
{1 \over \tilde{A}_d }
I_n(\gamma - d;d) & = &
\int_{R'} dx_0 \cdots dx_n dy_1 \cdots dy_n \,
\prod_{j=0}^n (1 + i x_j)^{-\gamma} (1 - i x_j)^{-\bar{\gamma}}
\nonumber \\
&& \times
\prod_{j=1}^n \prod_{l=0}^n | y_j - x_l|^{d-1}
\prod_{1 \le j < k \le n} |y_j - y_k|
\prod_{0 \le j < k \le n} |x_j - x_k|.
\end{eqnarray}

On the other hand, we know from \cite{Di05,An91} that
\begin{eqnarray}\label{66'}
&&
\int_{R'} dy_1 \cdots dy_n \,
\prod_{1 \le j < k \le n} |y_j - y_k|
\prod_{j=1}^n \prod_{l=0}^n | y_j - x_l|^{d - 1}
\nonumber \\
&& \qquad =
{(\Gamma(d))^{n+1} \over \Gamma((n+1)d) }
\prod_{0 \le j < k \le n} |x_j - x_k|^{2 d - 1},
\end{eqnarray}
so the r.h.s.~of (\ref{66}) is also equal to
\begin{equation}\label{67}
{(\Gamma(d))^{n+1} \over \Gamma((n+1)d) } I_{n+1}(\gamma;d).
\end{equation}
Equating (\ref{66}) and (\ref{67}) gives (\ref{cf}).
\hfill $\square$

\medskip
Iterating (\ref{cf}) with $I_0(\gamma,d) = 1$ reclaims the
gamma function evaluation \cite{Fo02}
\begin{equation}
n! I_n(\gamma;d) = 2^{d n(n-1) - 2 ({\rm Re} \, \gamma - 1) } \pi^n 
M_n(\bar{\gamma} - d(n-1)-1,\gamma -  d(n-1)-1,d)
\end{equation}
where $M_N(a,b,\lambda) $ is given by (\ref{4.16'}).

\subsection{A random three term recurrence}
Consider the rational function (\ref{cc2}). Suppose $\{q_j\}_{j=0,\dots,n}$
have the Dirichlet distribution (\ref{w1}) with equal parameters (\ref{al}), 
and suppose $c$ has the distribution (\ref{w2}). Suppose furthermore that
$\{y_j\}_{j=1,\dots,n}$ have distribution with measure
\begin{equation}\label{mh}
{1 \over I_n(\gamma - d;d) }
\prod_{j=1}^n (1 + i y_j)^{-\gamma + d} (1 - i y_j)^{-\bar{\gamma} + d}
\prod_{1 \le j < k \le n} | y_k - y_j|^{2d}.
\end{equation}
The marginal distribution of $\{x_j\}$ is then given by multiplying this 
with (\ref{w4}) and integrating $\{y_j\}$ over the region $R'$ (\ref{cc2'}). Using
(\ref{66'}) gives
\begin{equation}\label{red}
{1 \over I_{n+1}(\gamma;d) }
\prod_{j=0}^n (1 + i x_j)^{-\gamma} (1 - i x_j)^{-\bar{\gamma}}
\prod_{0 \le j < k \le n} |x_j - x_k|^{2d}.
\end{equation}
Hence with $p_{n+1}(z;\gamma;d)$ denoting the random monic polynomial of degree
$n+1$ with zeros at $\{x_j\}_{j=0,\dots,n}$ having distribution
(\ref{red}), we see that (\ref{cc2}) can be written
\begin{equation}\label{r10}
{p_{n+1}(z;\gamma;d) \over (z^2+1) p_n(z;\gamma - d; d) } =
{z - c \over q_0 (z^2 + 1) } - \sum_{j=1}^n {(q_j/q_0) \over z - y_j}.
\end{equation}

A companion identity to (\ref{r10}) is also required. For this purpose we introduce
the random rational function
\begin{equation}\label{2r2}
{\prod_{k=1}^{n-1} (z - u_k) \over \prod_{j=1}^n ( z - y_j) }
= \sum_{j=1}^n {\mu_j \over z - y_j}
\end{equation}
where $\{\mu_j\}$ have Dirichlet distribution
$$
{\Gamma(n d) \over (\Gamma(d))^n } \prod_{j=1}^n \mu_j^{d-1}.
$$
We know from the work of Dixon \cite{Di05} and Anderson \cite{An91} that the
conditional p.d.f.~of $\{u_k\}$ given $\{y_j\}$ is equal to
\begin{equation}\label{uk}
{\Gamma(n d) \over (\Gamma(d))^n }
{\prod_{1 \le j < k \le n - 1} (u_j - u_k) \over
\prod_{1 \le j < k \le n} (y_j - y_k)^{2 d - 1} }
\prod_{j=1}^{n-1} \prod_{k=1}^n | u_j - y_k|^{d - 1},
\end{equation}
provided
\begin{equation}\label{Rtt}
y_1 > u_1 > \cdots > y_{n-1} > u_{n-1} > y_n.
\end{equation}

It follows from (\ref{uk}) that if $\{y_j\}$ have distribution (\ref{mh}), then
the marginal distribution of $\{u_j\}$ is equal to
\begin{eqnarray*}
&& {\Gamma(n d) \over (\Gamma(d))^n} {1 \over I_n(\gamma - d; d) }
\prod_{1 \le j < k \le n - 1} (u_j - u_k)
\int_{\tilde{R}} dy_1 \cdots dy_n \,
\nonumber \\
&& \quad \times
\prod_{j=1}^n (1 + i y_j)^{-\gamma + d}
(1 - i y_j)^{-\bar{\gamma} + d} \prod_{j=1}^{n-1} \prod_{k=1}^n |u_j - y_k|^{d
-1}
\end{eqnarray*}
where $\tilde{R}$ is the region (\ref{Rtt}). According to (\ref{sgs}) this
can be evaluated as
$$
{1 \over I_{n-1}(\gamma -2d,d)}
\prod_{j=1}^{n-1}(1 + iu_j)^{-\gamma + 2d} (1 -iu_j)^{-\bar{\gamma} +2d}
\prod_{1 \le j < k \le n-1} |u_j - u_k|^{2 d}.
$$
We therefore conclude that (\ref{2r2}) can be written
\begin{equation}\label{r13}
{p_{n-1}(z;\gamma - 2d;d) \over p_n(z; \gamma - d;d) }
= \sum_{j=1}^n {\mu_j \over z - y_j}.
\end{equation}

Comparison of (\ref{r10}) and (\ref{r13}) 
implies $\{p_n(z;\gamma+(n-1)d;d)\}$ satisfy a
random three term recurrence.

\begin{thm}
With $B[\alpha,\beta]$ denoting the classical beta distribution, let
\begin{equation}\label{qoq}
b_n \sim B[ 2 {\rm Re} \, \gamma + nd - 1,nd] \quad (n \ne 0),
\qquad b_0=1,
\end{equation}
and let $c_n$ have the Cauchy distribution
\begin{equation}\label{qoq'}
{\Gamma(\gamma + n d) \Gamma(\bar{\gamma} + n d) \over \pi
2^{2(1 - nd - {\rm Re} \, \gamma)} \Gamma(2( {\rm Re} \, \gamma + nd) - 1)}
(1 + ic)^{-(\gamma + nd)} (1 - ic)^{-(\bar{\gamma} + nd)}  
\end{equation}
(this is (\ref{w2}) with $\gamma \mapsto \gamma + nd$). We have that
for $n=0,1,\dots,$
\begin{eqnarray}\label{qop}
\lefteqn{p_{n+1}(z;\gamma + nd;d)} \nonumber \\
&& \qquad = {(z - c_n) \over b_n} p_n(z;\gamma+(n-1)d;d) +
\Big ( 1 - {1 \over b_n} \Big ) (1 + z^2)
p_{n-1}(z;\gamma+(n-2)d;d),
\end{eqnarray}
where $p_0:=1$.
\end{thm}

\noindent
Proof. \quad In (\ref{r10}) and (\ref{r13}) replace $\gamma \mapsto \gamma + nd$.
We know that in general if $\{d_j\}_{j=0,\dots,n}$ have Dirichlet distribution
(\ref{w1}), then each $d_j$ has beta distribution 
B$[d_j, \sum_{l=0, l \ne j}^n
d_j]$. Using this fact it follows that in (\ref{r10}) we now have
$$
q_j \sim {\rm B}[d,(n-1)d] \: \: (j \ne 0), \qquad
q_0 \sim {\rm B}[2 {\rm Re} \, \gamma + nd - 1, nd],
$$
where in deriving the former use has also been made of (\ref{w3}), while in (\ref{r13}) 
$$
\mu_j \sim B[d,(n-1)d].
$$
The quantities are constrained by $\sum_{j=0}^n q_j = 1$, $\sum_{j=1}^n \mu_j = 1$.
Substituting (\ref{qop}) in (\ref{r10}) we thus see that (\ref{r13}) results, thereby
verifying the correctness of (\ref{qop}). 
\hfill $\square$

To relate this to the circular Jacobi $\beta$-ensemble (\ref{Mn}), we
note that with
$$
x_j = i {1 - e^{i \theta_j} \over 1 + e^{i \theta_j} }
\: \: (j \ne 0), \qquad x_0 = 0,
$$
the p.d.f.~(\ref{red}) with $\gamma \mapsto \gamma + 2d$ ($\gamma$ real)
extended to a measure via the multiplication by
$dx_1 \cdots dx_n$, becomes equal to (\ref{Mn}) with $a=2\gamma-2$ and
extended to a measure via the multiplication by $d \theta_1 \cdots
d \theta_n$. Thus the zeros of the polynomial $p_n(z;\gamma+(n-1)d;d)$,
with $\gamma$ real, $x_1,\dots,x_n$ say, under the mapping
\begin{equation}\label{us}
{x_j - i \over x_j + i } = e^{i \theta_j} \quad
(j=1,\dots,n)
\end{equation}
give for $\{\theta_j\}$ the distribution (\ref{Mn}) with 
$a=2 \gamma-2$.

As an illustration, let us consider the case $\gamma = 1$, $d=1$, which relates to
averaging over $U(N)$. There are a number of averages over $U(N)$ which are known
analytically. For example, with $p \in \mathbb Z_{>0}$, \cite{Dy62a}
$$
\langle | {\rm Tr} \, U^p|^2 \rangle_{U \in U(N)} =
\left \{ \begin{array}{ll} p, & 0 < p \le N \\
N, & p \ge N \end{array} \right.
$$
Since
$$
\langle | {\rm Tr} \, U^p |^2 \rangle_{U \in U(N)} =
\Big \langle \Big | \sum_{j=1}^N e^{ip \theta_j} \Big |^2 \Big \rangle_{U(N)}
$$
we can compute the Monte Carlo approximation
\begin{equation}\label{MC}
\langle | {\rm Tr} \, U^p|^2 \rangle_{U \in U(N)} =
{1 \over M} \sum_{k=1}^M \Big |
\sum_{j=1}^N \Big ( {x_j^{(k)} - i \over x_j^{(k)} + i } \Big )^p \Big |^2 +
O \Big ( {1 \over \sqrt{M}} \Big )
\end{equation}
where use has been made of (\ref{us}) and $x_j^{(k)}$ refers to the $j$th generation of
$p_N(z;N;1)$ from  (\ref{qop}). 

For the latter task, we read off from (\ref{qoq}) that
$$
b_n \sim {\rm B}[1+n,n].
$$
Also, by definition the Student $t$-distribution $T_\nu$ say, has p.d.f.~proportional to
$(1+t^2/\nu)^{-(\nu+1)/2}$ so
$$
c_n \sim {1 \over \sqrt{\nu} } T_\nu \Big |_{\nu = 2n+1}.
$$
Significantly, the zeros of the
lower order polynomials in the sequence $\{p_j(z;j;1)\}_{j=0,1,\dots,N}$
themselves allow us, via (\ref{us}), to sample from $U(j)$. Hence (\ref{MC}) can be calculated
for all values of $N$ less than the sought value within the same calculation.
Monte Carlo results obtained this way are presented in Table \ref{ta1}. The consistency
of these results is evident.

\begin{table}\label{ta1}
\begin{center}
\begin{tabular}{|c|c|c|c|c|c|}\hline
$p\backslash N$ & 2 & 3 & 4 & 5 \\\hline
1 & 0.98 & 0.99 & 0.99 & 0.98 \\\hline
2 & 2.02 & 2.00 & 2.05 & 1.99 \\ \hline
3 & 2.00 & 3.00 & 2.95 & 3.00 \\ \hline
4 & 1.96 & 3.00 & 3.97 & 4.01 \\ \hline
5 & 2.03 & 2.98 & 4.00 & 5.05 \\\hline
\end{tabular}
\end{center}
\caption{Computation of (\ref{MC}) with $M = 5,000$ and $p$ and $N$ as indicated}
\end{table}


\end{document}